\magnification=\magstep1
\overfullrule=0pt
\tolerance=700
\def\fnt#1{\csname !f-#1\endcsname}
\def\deffnt#1{\expandafter\font\csname !f-#1\endcsname}
\def\scal#1{ scaled \magstep#1}
\deffnt{rm}=cmr10
\deffnt{rm.1}=cmr10 \scal1
\deffnt{rm.2}=cmr10 \scal2
\deffnt{rm9}=cmr9
\deffnt{rm8}=cmr8
\deffnt{rm7}=cmr7
\deffnt{rm5}=cmr5
\deffnt{bf}=cmbx10
\deffnt{bf.1}=cmbx10 \scal1
\deffnt{bf.2}=cmbx10 \scal2
\deffnt{bf9}=cmbx9
\deffnt{bf7}=cmbx7
\deffnt{bf5}=cmbx5
\deffnt{it}=cmti10
\deffnt{it.1}=cmti10 \scal1
\deffnt{it9}=cmti9
\deffnt{it7}=cmti7
\deffnt{mi}=cmmi10
\deffnt{mi.1}=cmmi10 \scal1
\deffnt{mi9}=cmmi9
\deffnt{mi7}=cmmi7
\deffnt{mi5}=cmmi5
\deffnt{sy}=cmsy10
\deffnt{sy9}=cmsy9
\deffnt{ex}=cmex10
\deffnt{ex9}=cmex9
\deffnt{sl}=cmsl10
\deffnt{sl8}=cmsl8
\deffnt{sf}=cmss10
\deffnt{sf.2}=cmss10 \scal2
\deffnt{sc}=cmcsc10
\deffnt{sc9}=cmcsc9
\deffnt{sc.1}=cmcsc10 \scal1
\deffnt{mb}=msbm10
\deffnt{mb7}=msbm7
\deffnt{mb5}=msbm5
\newfam\mbfam
\textfont\mbfam=\fnt{mb}
\scriptfont\mbfam=\fnt{mb7}
\scriptscriptfont\mbfam=\fnt{mb5}
\def\mb{\fam\mbfam\fnt{mb}}
\newfam\lgfam
\textfont\lgfam=\fnt{mi.1}
\scriptfont\lgfam=\fnt{mi9}
\scriptscriptfont\lgfam=\fnt{mi7}

\newfam\smfam
\textfont\smfam=\fnt{mi9}
\scriptfont\smfam=\fnt{mi7}
\scriptscriptfont\smfam=\fnt{mi5}
\deffnt{gt}=eufm10
\deffnt{gt5}=eufm5

\def\cmd#1{\csname #1\endcsname}
\def\defcmd#1{\expandafter\def\csname #1\endcsname}
\def\edefcmd#1{\expandafter\edef\csname #1\endcsname}
\def\ifundefined#1{\expandafter\ifx\csname #1\endcsname\relax}
\edef\warning#1{{\newlinechar=`|\message{|*********** #1}}}
\def\npar{\relax\par\noindent}
\newbox\xbox
\def\ifhempty#1#2#3{\setbox\xbox=\hbox{#1}\ifdim\wd\xbox=0pt\relax
                                                           #2\else #3\fi}
\def\lph#1{\hbox to 0pt{\hss #1}}
\def\rph#1{\hbox to 0pt{#1\hss}}
\def\ulph#1{\vbox to 0pt{\vss\hbox to 0pt{\hss #1}}}
\def\urph#1{\vbox to 0pt{\vss\hbox to 0pt{#1\hss}}}
\def\dlph#1{\vbox to 0pt{\hbox to 0pt{\hss #1}\vss}}
\def\drph#1{\vbox to 0pt{\hbox to 0pt{#1\hss}\vss}}
\def\nasp{\spacefactor=1000}
\def\small{\def\rm{\fnt{rm9}\fam\smfam}\def\it{\fnt{it9}}\def\bf{\fnt{bf9}}%
\textfont1=\fnt{mi9}\textfont2=\fnt{sy9}\textfont3=\fnt{ex9}%
\baselineskip=.8\baselineskip\parindent=.8\parindent%
\abovedisplayskip=.6\abovedisplayskip\belowdisplayskip=.6\belowdisplayskip%
\stateskip=.6\stateskip
\rm}
\def\nohyphen{\hyphenpenalty=5000 \exhyphenpenalty=5000
 \tolerance=5000 \pretolerance=5000}
\def\center{\leftskip=0pt plus 1fill\relax \rightskip=\leftskip\relax
 \parfillskip=0pt\relax \parindent=0pt\relax \nohyphen}
\def\\{\ifhmode\hfil\break\else\ifvmode\vskip\baselineskip\fi\fi}
\def\-{\hfill}
\def\vsp#1{\vadjust{\kern #1}}

\def\vbreak#1#2{\par\ifdim\lastskip<#2\removelastskip\penalty#1\vskip#2\fi}
\def\vbreakn#1#2{\vbreak{#1}{#2}\npar}

\newcount\multilabels \newcount\undeflabels
\newread\filelabelsold \newwrite\filelabels
\edef\filelabelsname{\jobname.lbl}
\def\getlabels{\openin\filelabelsold=\filelabelsname\relax
 \ifeof\filelabelsold 
  \warning{The file of labels \filelabelsname\ does not exist}
  \else\closein\filelabelsold
   {\catcode`/=0 \globaldefs=1 \input \filelabelsname\relax}\fi}
\def\openlabels{\immediate\openout\filelabels=\filelabelsname\relax}
\def\closelabels{\closeout\filelabels}

\def\label#1#2#3{
\ifx\printlabel Y\lnote{\fnt{rm7}#1}\fi%
\immediate\write\filelabels{%
/defcmd{!tl-#1}{#2}/defcmd{!nl-#1}{#3}}%
\ifundefined{!dl-#1}\defcmd{!dl-#1}{!earlier-defined!}%
\else\warning{Label #1 is multiply defined.}\advance\multilabels by 1\fi%
\edefcmd{!tl-#1}{#2}\edefcmd{!nl-#1}{#3}}
\def\xlabel#1#2#3#4{
\ifx\printlabel Y\ldnote{\fnt{rm7}#1}\fi%
\immediate\write\filelabels{%
/defcmd{!tl-#1}{#2}/defcmd{!nl-#1}{#3}/defcmd{!xl-#1}{#4}}%
\ifundefined{!dl-#1}\defcmd{!dl-#1}{!earlier-defined!}%
\else \warning{LABEL #1 is multiply defined.} \advance\multilabels by 1\fi%
\edefcmd{!tl-#1}{#2}\edefcmd{!nl-#1}{#3}\edefcmd{!xl-#1}{#4}}
\def\undelabel#1{\warning{Label #1 has not been defined.}%
\global\advance\undeflabels by 1{?#1?}}
\def\rfrn#1{\ifundefined{!nl-#1}\undelabel{#1}\else\cmd{!nl-#1}\fi}

\def\rfrx#1{\ifundefined{!xl-#1}\undelabel{#1}%
 \else\cmd{!xl-#1}\fi}
\def\rfr#1{\ifundefined{!tl-#1}\undelabel{#1}%
 \else\cmd{!tl-#1}~\cmd{!nl-#1}\fi}

\def\title#1{\vbreak{-200}{1cm}{\center\fnt{bf.1}#1\par}}
\def\author#1{\vbreak{100}{5mm}{\center\fnt{rm}#1\par}}
\def\titledate{\vbreak{100}{5mm}\centerline{\fnt{it}\date}}

\newskip\absskip
\absskip=20mm
\long\def\abstract#1{\vbreak{50}{5mm}\centerline{\fnt{bf}Abstract}
 \kern3mm{\leftskip=\absskip\rightskip=\leftskip
 \small #1\vbreak{-9000}{4mm}}}

\def\support#1{\footnote{}{{\small#1\hfil}}}

\ifx\headdate Y{\headline={{\fnt{rm7}\date}\hfil}}\fi
\def\date{\ifcase\month%
\or January\or February\or March\or April\or May\or June\or July%
\or August\or September\or October\or November\or December\fi%
~\number\day,~\number\year}

\def\myaddress{\rm%
\setbox\xbox=\hbox{Department of mathematics}%
\vtop{\hsize=\wd\xbox\parindent=0pt
Department of Mathematics\\
The Ohio State University\\
Columbus, OH 43210, USA\\
\hbox to \wd\xbox{{\it e-mail}: leibman@math.ohio-state.edu\hss}}}


\newcount\fnotenum
\def\fnote#1{\global\advance\fnotenum by 1
 {\parindent=10mm \parfillskip=0pt plus 1fill\relax
  \footnote{$^{(\the\fnotenum)}$}{{\small #1}}}}

\def\lnote#1{\ifvmode\ulph{#1 }\else\vadjust{\ulph{#1 }}\fi}
\def\ldnote#1{\ifvmode\ulph{#1 }\else\vadjust{\dlph{#1 }}\fi}
\newcount\sectionnum \sectionnum=-1 
\newcount\subsectionnum \subsectionnum=0 
\newcount\subsubsectionnum \subsubsectionnum=0
\newskip\sectionskip \sectionskip=6mm
\newskip\aftersectionskip \aftersectionskip=3mm
\newcount\sectionpenalty \sectionpenalty=-9000
\def\xsection#1{\global\advance\sectionnum by 1%
\global\subsectionnum=0\global\subsubsectionnum=0%
\ifx\nonewsecequation Y\else\global\equationnum=0\fi%
\ifx\nonewsecstate Y\else\global\statenum=0\fi%
\label{#1}{Section}{\the\sectionnum}}
\def\section#1#2{%
\vbreak{\sectionpenalty}{\sectionskip}\xsection{#1}%
\vbox{\ifx\nosectioncenter Y\npar\else\center\fi\fnt{bf}
\ifx\nosectionnum Y\else\rfrn{#1}.{\nasp} \fi#2\par}%
\vbreak{10000}{\aftersectionskip}}
\newskip\subsectionskip \subsectionskip=3mm
\newcount\subsectionpenalty \subsectionpenalty=-7000
\def\thesubsectionnum{\ifx\nosection Y\the\subsectionnum
\else\the\sectionnum.\the\subsectionnum\fi}
\def\xsubsection#1{\global\advance\subsectionnum by 1%
\global\subsubsectionnum=0%
\label{#1}{subsection}{\thesubsectionnum}}
\def\subsection#1{%
\vbreak{\subsectionpenalty}{\subsectionskip}
\ifx\nosubsectionnum Y\par\else
\noindent\xsubsection{#1}{\fnt{bf}\rfrn{#1}.{\nasp}}\fi}
\def\thesubsubsectionnum{\thesubsectionnum.\the\subsubsectionnum}
\def\xsubsubsection#1{\global\advance\subsubsectionnum by 1%
\label{#1}{subsection}{\thesubsubsectionnum}}
\def\subsubsection#1{%
\vbreak{-5000}{1.5mm}\noindent\xsubsubsection{#1}{\fnt{bf}\rfrn{#1}.{\nasp}}}
\newcount\equationnum   
\def\theequation{\ifx\nosection Y\the\equationnum
\else\the\sectionnum.\the\equationnum\fi}
\def\equn#1#2{\global\advance\equationnum by 1%
\label{#1}{Equation}{\theequation}%
$$\vcenter{\equalign{#2}}\eqno(\theequation)$$}
\def\equ#1{$$\vcenter{\equalign{#1}}$$}
\def\equalign#1{\let\\=\cr\let\-=\hfill
\ialign{&\hfil$\dsp ##$\hfil\cr#1\crcr}}

\def\lequalign#1{\let\\=\cr\let\-=\hfill
\ialign{&\hbox to \hsize{$\dsp ##$\hfil}\cr#1\crcr}}
\def\matalign#1#2{{\let\\=\cr\let\-=\hfill%
\ialign{\hfil$##$\hfil&&#1\hfil$##$\hfil\cr#2\crcr}}}

\def\frfr#1{\hbox{(\rfrn{#1})}}
\newskip\stateskip\stateskip=2mm
\newcount\statenum\statenum=0
\def\thestate{\ifx\afterstate Y\the\statenum\else\thesubsectionnum\fi}
\def\state#1#2#3#4{\global\advance\statenum by 1\xlabel{#2}{#1}{\thestate}{#3}%
{\bf #1\ifhempty{#3}{\ifx\afterstate Y\ \thestate\fi}{ \box\xbox}.}{\nasp} #4%
\begingroup}
\def\nstate#1#2{{\bf #1\ifhempty{#2}{}{ \box\xbox}.}{\nasp}\begingroup}
\def\endstate{\par\endgroup\vbreak{-7000}{\stateskip}}
\def\ltheorem#1#2#3{\state{Theorem}{#1}{#2}{#3}\it}
\def\theorem#1#2#3{\vbreakn{-7000}{\stateskip}\ltheorem{#1}{#2}{#3}}
\def\endtheorem{\endstate}
\def\llemma#1#2#3{\state{Lemma}{#1}{#2}{#3}\it}
\def\lemma#1#2#3{\vbreakn{-7000}{\stateskip}\llemma{#1}{#2}{#3}}
\def\endlemma{\endstate}
\def\lproposition#1#2#3{\state{Proposition}{#1}{#2}{#3}\it}
\def\proposition#1#2#3{\vbreakn{-7000}{\stateskip}\lproposition{#1}{#2}{#3}}
\def\endproposition{\endstate}
\def\lcorollary#1#2#3{\state{Corollary}{#1}{#2}{#3}\it}
\def\corollary#1#2#3{\vbreakn{-7000}{\stateskip}\lcorollary{#1}{#2}{#3}}
\def\endcorollary{\endstate}


\def\lremark#1{\nstate{Remark}{#1}}
\def\remark#1{\vbreakn{-7000}{\stateskip}\lremark{#1}}
\def\endremark{\endstate}






\def\lproof#1{\nstate{Proof}{#1}}
\def\proof#1{\vbreakn{-7000}{\stateskip}\lproof{#1}}

\def\endproof{\endprr\endstate}
\def\enprule{\vrule height1mm depth1mm width2mm}
\def\endprr{\discretionary{}{\kern\hsize}{\kern 3ex}\llap{\enprule}}
\def\frgdsp{\par\penalty10000
 \vskip-\belowdisplayskip\kern-2mm\noindent\hbox to \hsize{\hfil}}
\newcount\biblnum \newcount\biblpenalty
\def\bibliography#1{\vbreak{-9000}{6mm}
\biblpenalty=10000\biblnum=0%
\line{\fnt{bf}\ifhempty{#1}{Bibliography}{\box\xbox}\hfil}\par\kern2mm
\begingroup\parskip=0pt\parindent=0pt\frenchspacing}
\def\endbibliography{\par\endgroup}
\newskip\biblleft \biblleft=10mm
\def\bibxitem#1/#2(#3) #4{\advance\biblnum by 1%
\vbreak{\biblpenalty}{.8mm}\biblpenalty=-9000%
\hangindent=\biblleft \hangafter=1
\noindent{\small\rlap{\brfr{#1}}\hskip\hangindent
\xlabel{#1}{#3}{\the\biblnum}{#2}#4}}
\def\bibitem#1/#2(#3) #4#5{\bibxitem#1/#2(#3) {{\frenchspacing #4}, #5}}
\def\bibx#1/#2 {\advance\biblnum by 1{#1/#2}
\xlabel{#1}{text}{\the\biblnum}{#2}}
\def\bibart#1/#2 a:#3 t:#4 j:#5 n:#6 y:#7 p:#8 *{%
\bibitem{#1}/{#2}(paper) {#3}{{#4}, {\it #5\/}
   {\bf #6} (#7), #8.}}
\def\bibartp#1/#2 a:#3 t:#4 *{%
\bibitem{#1}/{#2}(paper) {#3}{{#4}, {in preparation}.}}
\def\bibarts#1/#2 a:#3 t:#4 *{%
\bibitem{#1}/{#2}(paper) {#3}{{#4}, {submitted}.}}
\def\bibarta#1/#2 a:#3 t:#4 j:#5 *{%
\bibitem{#1}/{#2}(paper) {#3}{{#4}, {\rm to appear in \it #5}.}}
\def\bibartn#1/#2 a:#3 t:#4 j:#5 *{%
\bibitem{#1}/{#2}(paper) {#3}{{#4}, #5.}}
\def\bibartpr#1/#2 a:#3 t:#4 *{%
\bibitem{#1}/{#2}(paper) {#3}{{#4}, {preprint}.}}
\def\bibartpc#1/#2 a:#3 *{%
\bibitem{#1}/{#2}(paper) {#3}{private communications}.}
\def\bibartx#1/#2 a:#3 t:#4 x:#5 *{%
\bibitem{#1}/{#2}(paper) {#3}{{#4}, #5.}}
\def\bibook#1/#2 a:#3 t:#4 i:#5 *{%
\bibitem{#1}/{#2}(book) {#3}{{\it #4\/}, #5.}}
\def\no#1{no.{\nasp}~{#1}}

\def\brfr#1{\hbox{\rm[\ifx\biblNUM Y\rfrn{#1}\else\rfrx{#1}\fi]}}
\def\start{
 \warning{|****** Start ******|}
 \getlabels\openlabels\multilabels=0\undeflabels=0
 \null}
\def\finish{\par\closelabels
 \warning{|****** Finish ******}
 \ifnum\multilabels>0\warning{\the\multilabels_multidefined labels!}\fi
 \ifnum\undeflabels>0\warning{\the\undeflabels_undefined labels!}
                     \warning{*************************}
                     \warning{Try to run TeX once again!}
                     \warning{*************************}\fi}
\output={\shipout\vbox{\makeheadline\pagebody\makefootline}
\ifx\doublepage Y\shipout\vbox{}\fi
 \advancepageno
 \ifnum\outputpenalty>-20000\else\dosupereject\fi}
\def\sdup#1#2{{\scr #1 \atop \scr #2}}

\def\Gal{\lnote{{\bf V\kern 1cm\relax}}}
\def\comment#1{\ifhmode\\\fi
\hbox to \hsize{\hss\vbox{\advance\hsize by 10mm
\baselineskip=.7\baselineskip\npar\ulph{\fnt{bf.2}V}\fnt{bf7}#1}}}

\def\rest#1{\raise-2pt\hbox{$|_{#1}$}}

\def\notdvd{\mathrel{\setbox\xbox=\hbox{$\big|$}%
\hbox to 0pt{\hbox to\wd\xbox{\hss/\hss}\hss}\big|}}

\def\frac#1#2{{#1\over #2}}
\def\matr#1{{\let\\=\cr\left(\matrix{#1\crcr}\right)}}
\def\vect#1{{\let\\=\cr\left(\matrix{#1\crcr}\right)}}
\def\smatr#1{{\baselineskip=2pt\lineskip=2pt
\left(\vcenter{\let\\=\cr\let\-=\hfill
\ialign{\hfil$\scr##$\hfil&&\hfil\kern2pt$\scr##$\hfil\cr#1\crcr}}\right)}}

\def\svd{\vbox to 2.4mm{\nolineskip
\kern.3pt\hbox{.}\vfil\hbox{.}\vfil\hbox{.}\kern.3pt}}
\def\vd{\vbox to 3.2mm{\nolineskip
\kern1pt\hbox{.}\vfil\hbox{.}\vfil\hbox{.}\kern1pt}}
\def\lvd{\vbox to 7mm{\nolineskip
\kern1pt\hbox{.}\vfil\hbox{.}\vfil\hbox{.}\vfil\hbox{.}
\vfil\hbox{.}\vfil\hbox{.}\kern1pt}}
\def\rvd{\vbox to 2.4mm{\nolineskip
\kern.3pt\hbox{.}\vfil\hbox{\kern3.5pt.}\vfil\hbox{\kern7pt.}\kern.3pt}}
\def\lld{\hbox to 5mm{\kern2pt.\hfil.\hfil.\kern2pt}}

\def\comp{\mathord{\hbox{$\scr\circ$}}}
\def\semprod{\mathrel{\times\kern-2pt
\vbox{\hrule width.5pt height4pt depth0pt\kern.5pt}\kern1pt}}

\def\dsc{\discretionary{}{}{}}
\long\def\omit#1\endomit{\par\vbox{%
\hrule\vskip2mm\hfil\vdots\hfil\vskip2mm\hrule}}
\long\def\ignore#1\endignore{}
\def\R{{\mb R}}

\def\Z{{\mb Z}}
\def\N{{\mb N}}

\def\PP{\raise2.2pt\hbox{\fnt{mi.1}\char"7D}}

\def\Id{\mathop{\hbox{\rm Id}}}


\def\nolineskip{\baselineskip=0pt\lineskip=0pt}

\let\dsp=\displaystyle
\let\txt=\textstyle

\let\scr=\scriptstyle

\let\sln=\subset
\let\sle=\subseteq

\let\sm=\setminus

\let\col=\colon

\let\ld=\ldots
\let\vd=\vdots
\let\cd=\cdot
\let\cds=\cdots
\let\pus=\emptyset
\let\ra=\longrightarrow

\let\alf=\alpha
\let\bet=\beta
\let\Gam=\Gamma
\let\gam=\gamma
\let\del=\delta

\let\phi=\varphi

\let\eps=\varepsilon

\let\ro=\rho


\let\printlabel=N
\let\afterstate=Y
\absskip=0mm
\subsectionskip=5mm

\def\thestate{\the\sectionnum.\the\statenum}

\def\inr#1{[1,#1]}
\def\lip{\|}
\def\rip{\|}
\def\cF{{\cal F}}

\def\TT{{\mb T}}

\def\dist{\mathop{\hbox{\rm dist}}}
\def\ndeg{\mathop{\hbox{\rm n-deg}}}
\def\Pol{\mathord{\hbox{\rm Pol}}}
\def\Polz{\mathord{\hbox{\rm Pol}}^{0}}

\def\GPolzo{\mathord{\hbox{\rm GPol}}^{0}}
\def\Lin{\mathord{\hbox{\rm Lin}}}
\def\rsub#1{\mathclose{\vtop{\kern-4pt\hbox{$\downarrow_{#1}$}}}}
\def\IP{IP}
\def\IPz{IP$^{*}$}
\def\IPr#1{IP$_{#1}$}
\def\IPzr#1{IP$_{#1}^{*}$}

\def\IPzn{IP$_{0}^{*}$}
\def\VIP{VIP}
\def\VIPz{VIP$^{*}$}
\def\VIPr#1{VIP$_{#1}$}
\def\VIPzr#1{VIP$_{#1}^{*}$}

\def\VIPzn{VIP$_{0}^{*}$}
\def\hPhi{\widehat{\Phi}}
\def\tPhi{\widetilde{\Phi}}
\def\hphi{\widehat{\phi}}
\def\RP#1{\hbox{\bf RP$^{[#1]}$}}
\let\Del=\Delta
\def\Delz{$\Del^{*}$}
\def\Delzn{$\Del^{*}_{0}$}
\def\sbc#1{disjoint \ifhempty{#1}{}{#1-}subcollection}
\def\cbc{q}
\def\smp{t}
\def\PP#1{#1}
\start
\title{IP$_{\hbox{\bf r}}^{\hbox{\bf *}}$-recurrence and nilsystems}
\author{V.~Bergelson and A.~Leibman}
\support{Partially supported by NSF grants DMS-1162073 and DMS-1500575.}
\titledate
\abstract{By a result due to Furstenberg, 
a homeomorphism $T$ of a compact space is distal if and only if it possesses the property of \IPz-recurrence,
meaning that for any $x_{0}\in X$, for any open neighborhood $U$ of $x_{0}$, 
and for any sequence $(n_{i})$ in $\Z$, 
the set $R_{U}(x_0) =\{n\in\Z:T^{n}x_{0}\in U\}$ has
a non-trivial intersection with the set of finite sums 
$\{n_{i_{1}}+n_{i_{2}}+\cds+n_{i_{s}}: i_{1}<i_{2}<\ld<i_{s},\ s\in\N\}$. 
We show that translations on compact nilmanifolds (which are known to be distal) 
are characterized by a stronger property of \IPzr{r}-recurrence,
which asserts that for any $x_{0}\in X$ and any neighborhood $U$ of $x_{0}$ 
there exists $r\in\N$ such that for any $r$-element sequence $n_{1},\ld,n_{r}$ in $\Z$
the set $R_{U}(x_{0})$ has a non-trivial intersection with the set 
$\{n_{i_{1}}+n_{i_{2}}+\cds+n_{i_{s}}: i_{1}<i_{2}<\ld<i_{s},\ s\leq r\}$.
We also show that the property of \IPzr{r}-recurrence is equivalent to an ostensibly much stronger
property of polynomial \IPzr{r}-recurrence. 
(This should be juxtaposed with the fact that for general distal transformations 
the polynomial \IPz-recurrence is strictly stronger than the \IPz-recurrence.)}
\section{S-Int}{Introduction}
Let $(X,T)$ be a topological dynamical system,
meaning that $X$ is a compact metric space and $T$ is a self-homeomorphism of $X$.
Given a point $x_{0}\in X$ and an open neighborhood $U$ of $x_{0}$,
define $R_{U}(x_{0})=\bigl\{n\in\Z:T^{n}x_{0}\in U\bigr\}$,
{\it the set of returns\/} of $x_{0}$ into $U$.
Sets of returns reflect the properties of topological system,
and it is of interest to characterize (and/or distinguish between) dynamical systems 
by arithmetic properties of these sets.
An example of this kind is provided by a theorem of Furstenberg
on sets of returns in distal systems.
A system $(X,T)$ is said to be {\it distal\/} 
if for any distinct $x,y\in X$, $\inf_{n\in\Z}\dist(T^{n}x,T^{n}y)>0$.
Given a sequence $n_{1},n_{2},\ld$ in $\Z$,
the set $\bigl\{n_{i_{1}}+\cds+n_{i_{s}}:s\in\N,\ i_{1}<\cds<i_{s}\bigr\}$ 
of finite sums of distinct elements of this sequence is called an {\it \IP-set}.
A subset $E$ of $\Z$ is called an {\it \IPz-set} if it intersects every \IP-set.
Furstenberg's theorem says that 
distal systems are characterized by {\it the \IPz-recurrence property}:
\theorem{P-IP}{}{(\brfr{F-book}, Theorem 9.11)}
A system $(X,T)$ is distal if and only if for any $x_{0}\in X$ and any open neighborhood $U$ of $x_{0}$
the set of returns $R_{U}(x_{0})$ is an \IPz-set.
\endtheorem

Another relevant example involves translations on compact abelian groups.
{\it A set of differences\/} is a set of the form $\bigl\{n_{i}-n_{j},\ j<i\bigr\}$,
where $(n_{i})$ is an infinite sequence in $\Z$;
a subset $E$ of $\Z$ is said to be {\it a \Delz-set\/}
if it has a nonempty intersection with every set of differences in $\Z$.
A point $x$ in a system $(X,T)$ is said to be {\it almost automorphic}
if for any sequence $(n_{i})$ in $\Z$,
$T^{n_{i}}x\ra y$ implies $T^{-n_{i}}y\ra x$.
It is shown in \brfr{F-book}, Theorem~9.13, 
that a system has the \Delz-recurrence property
(that is, that every set of returns in the system is a \Delz-set)
if and only if every point in the system is almost automorphic.
Next, by a theorem of Veech 
(see \brfr{Veech}, Theorem~1.2; see also \brfr{AGN})
every point of a minimal%
\fnote{A system $(X,T)$ is {\it minimal\/} if it has no proper closed subsystems,
or, equivalently, if the orbit of every point of $X$ is dense in $X$.}
system $(X,T)$ is almost automorphic
if and only if the family $\{T^{n},\ n\in\Z\}$ of powers of $T$ is equicontinuous.
Now, it is not hard to see that for a minimal $T$
the family $\{T^{n},\ n\in\Z\}$ is equicontinuous
if and only if $(X,T)$ is isomorphic to a translation on a compact abelian group%
\fnote{The ``only if'' implication follows from the fact that for any $x_{0}\in X$ 
one can define an additive group structure on the orbit $\{T^{n}x_{0},\ n\in\Z\}$ 
by $T^{n}x_{0}+T^{m}x_{0}=T^{n+m}x_{0}$, $n,m\in\Z$, 
and then extend it, with the help of equicontinuity, to all of $X$. 
This makes $X$ a compact abelian group on which $T$ acts as a minimal translation.}.
Thus, the recurrence property characterizing minimal group translations is that of \Delz.

Our goal in this paper is to provide a similar characterization of {\it nilsystems\/},
namely, systems of the form $(X,T)$ 
where $X$ is a {\it nilmanifold\/} 
(a compact homogeneous space of a nilpotent Lie group $G$)
and $T$ is a {\it niltranslation} 
(a translation on $X$ defined by an element of $G$).
The motivation for this study comes from 
the fact that nilsystems are intrincically related to various problems
arising in ergodic theory of multiple recurrence, combinatorics, and number theory,
and understanding the recurrence properties of niltranslations 
leads to interesting applications in these areas.
It is well known that nilsystems are distal (see \brfr{AGH}, \brfr{Keyns1}, \brfr{Keyns2}),
and thus are \IPz-recurrent;
however, not every distal system is a nilsystem,
and thus there must be a stronger than \IPz\ property of recurrence that characterizes them.

For an integer $r\in\N$ and an $r$-element sequence $n_{1},\ld,n_{r}$ in $\Z$,
we call the set
$\bigl\{n_{i_{1}}+\cds+n_{i_{s}}:1\leq s\leq r,\ i_{1}<\cds<i_{s}\bigr\}$
of sums of distinct elements of this sequence an {\it \IPr{r}-set}.
A set $E\sle\Z$ is called an {\it \IPzr{r}-set}
if it has a nonempty intersection with every \IPr{r}-set in $\Z$.
We say that a set is an {\it \IPzn-set\/} if it is an \IPzr{r}-set for some $r\in\N$.
\IPzn-sets form a proper subfamily of the family of \IPz-sets:
clearly, every \IPzn-set is \IPz, but not vice~versa%
\fnote{To see this, it is enough to exhibit an \IPz-set $S$ which is not an IP-set.
One can take, for example $S=\bigcup_{r=1}^{\infty}S_{r}$, 
where $S_{r}=\{2^{2^{r}}, 2\cd2^{2^{r}},3\cd2^{2^{r}},\ld,r\cd2^{2^{r}}\}$, $r\in\N$.  
Since for each $r$, $S_{r}$ is a dilation of the set $\{1,2,...,r\}$, 
$S$ contains arbitrarily large \IPr{r}-sets,
but it contains no IP-sets
since the distances between consecutive elements of $S$ form a non-decreasing sequence which tends to infinity.}.
A special class of nilsystems is provided by affine skew product transofrmations of tori%
\fnote{An affine skew product transformation of the $k$-dimensional torus $\TT^{k}=\R^{k}/\Z^{k}$ 
is defined by the formula $T(x_{1},\ld,x_{k})=(x_{1}+\alf_{1},x_{2}+a_{2,1}x_{1}+\alf_{2},\ld,
x_{k}+a_{k,k-1}x_{k-1}+\cds+a_{k,1}x_{1}+\alf_{k})$
with $\alf_{i}\in\TT$ and $a_{i,j}\in\Z$.};
it follows from \brfr{B-ultra}, Theorem~7.7, 
that every such system has the \IPzn-recurrence property:
for every $x_{0}\in\TT^{k}$ and any open neighborhood $U$ of $x_{0}$ 
the set of returns $R_{U}(x_{0})$ is an \IPzn-set.
On the other hand, 
one can show that not every distal system is \IPzn-recurrent (see \brfr{bic}, Section~1).
It is tempting to conjecture that it is the \IPzn-recurrence property that characterizes the nilsystems.
This, however, cannot be exactly so:
any recurrence property must be stable under passing to inverse limits
whereas inverse limits of nilsystems do not have to be nilsystems.
Let us define {\it a pre-nilsystem\/} as the inverse limit of a sequence of nilsystems.
(Notice that, in contrast with the definition of the so-called {\it pro-nilsystems},
in the definition of pre-nilsystems
we don't require the nilpotency class of the nilsystems in the sequence to be bounded.)
The following result provides a characterization of pre-nilsystems
in terms of \IPzn-recurrence:
\theorem{P-IPn}{}{}
Any pre-nilsystem (and so, any nilsystem) is \IPzn-recurrent.
Any \IPzn-recurrent system is a disjoint union of pre-nilsystems.
\endtheorem

\remark{}
In analogy with \IPzn-sets,
one can define \Delzn-sets as those having a nonempty intersection
with every large enough finite set of differences.
In contrast with \IPz/\IPzn-recurrence,
the classes of \Delz- and \Delzn-recurrent systems coincide.
(These are translations of compact abelian groups.)
\endremark

The second statement of \rfr{P-IPn} is an easy corollary of the results from \brfr{HKM}.
To prove the first statement, we use a coordinate approach.
On any nilmanifold $X$ one has natural coordinates
such that under the action of a niltranslation $T$
the sequence of coordinates of the image $T^{n}x_{0}$ of any point $x_{0}\in X$
is given by {\it generalized polynomials\/} (see \brfr{sko}, Theorem~A).
We therefore need to deal with images of \IP-sets under generalized polynomial mappings;
these images form a subclass of {\it generalized polynomial \IP-sets}.
Conventional \IP- and \IPr{r}-sets in $\Z$
can be viewed as the images of mappings $\phi\col\cF(A)\ra\Z$
from the semigroup $\cF(A)$ of finite subsets of $A$, for $A=\N$ and, respectively, for $A=\{1,\ld,r\}$,
defined by $\phi(\alf)=\sum_{i\in\alf}a_{i}$.
Such a mapping $\phi$ is ``linear'' in the following sense:
$\phi(\alf\cup\bet)=\phi(\alf)+\phi(\bet)$ whenever $\alf,\bet\in\cF(A)$ are disjoint.
Let $H$ be an additive abelian group;
one can introduce the notion of polynomial mappings $\cF(A)\ra H$ as follows.
For a mapping $\phi\col\cF(A)\ra H$ and a set $\bet\in\cF(A)$
let {\it the $\bet$-derivative\/} $D_{\bet}\phi$ be the mapping $\cF(A\sm\bet)\ra H$
defined by $D_{\bet}\phi(\alf)=\phi(\alf+\bet)-\phi(\alf)$.
Then we say that a mapping $\phi\col\cF(\{1,\ld,r\})\ra H$ is {\it polynomial of degree $\leq d$} 
if for any disjoint $\bet_{0},\bet_{1},\ld,\bet_{d}\in\cF(\{1,\ld,r\})$,
$D_{\bet_{0}}D_{\bet_{1}}\cds D_{\bet_{d}}\phi=0$.
(See \brfr{phj}, Section~8.1.)
Examples of quadratic (that is, of degree $\leq 2$) polynomial mappings are,
in increasing generality, 
$\phi(\alf)=\bigl(\sum_{i\in\alf}a_{i}\bigr)^{2}$,
$\phi(\alf)=\bigl(\sum_{i\in\alf}a_{i}\bigr)\bigl(\sum_{i\in\alf}b_{i}\bigr)
=\sum_{i,j\in\alf}a_{i}b_{j}$,
and $\phi(\alf)=\sum_{i,j\in\alf}c_{i,j}$,
where $a_{i},b_{j},c_{i,j}\in H$.
{\it Generalized polynomial mappings\/} are the mappings
built from (conventional) polynomial mappings 
using the operations of addition, multiplication, and taking the integer part.
(An example is $\phi=\bigl[[\phi_{1}]\phi_{2}+\phi_{3}\bigr]\phi_{4}+[\phi_{5}][\phi_{6}]\phi_{7}$,
which is comprised of the polynomial mappings $\phi_{1},\ld,\phi_{7}$.)
Let us say that a generalized polynomial mapping $\phi$  has {\it total degree\/} $\leq D$
if the sum $\sum_{i}\deg\phi_{i}$ of the degrees 
of all the ``conventional'' polynomial mappings $\phi_{i}$ of which $\phi$ is comprised
does not exceed $D$,
and let us say that a generalized polynomial mapping is {\it constant free\/} 
if all the $\phi_{i}$ vanish at $\pus$: $\phi_{i}(\pus)=0$.
Let us also say that a generalized polynomial mapping is {\it open\/}
if it is contained in the ideal generated by the conventional constant-free polynomials
of the ring of constant-free generalized polynomials.
(In other words, a generalized polynomial is open 
if it contains no ``closed'' summands of the form $[\phi_{1}]\cds[\phi_{k}]$,
where $\phi_{i}$ are generalized polynomials.)
For $x\in\R$ let $\lip x\rip=\dist(x,\Z)$.
The following result of Diophantine nature,
which we use to prove \rfr{P-IPn}, 
is of independent interest:
\theorem{P-Igpssub}{}{(Cf.\nasp\ \rfr{P-gpssub} below.)}
For any $D\in\N$ and $\eps>0$ there exists $r=r(D,\eps)\in\N$ 
such that for any open constant-free generalized polynomial mapping $\phi\col\{1,\ld,r\}\ra\R$
of total degree $\leq D$
there exists a nonempty $\alf\sle\{1,\ld,r\}$ for which $\lip\phi(\alf)\rip<\eps$.
\endtheorem

The {\it\VIP-sets\/} in $\Z^{l}$
are defined as the images $\bigl\{\phi(\alf):\alf\in\cF(\N),\ \alf\neq\pus\bigr\}$ 
of polynomial mappings $\phi\col\cF(\N)\ra\Z^{l}$ with $\phi(\pus)=0$,
and we say that a set $E\sle\Z^{l}$ is {\it a \VIPz-set\/}
if $E$ has a nonempty intersection with every \VIP-set in $\Z^{l}$.
Similarly, for all $d,r\in\N$, we define {\it \VIPr{d,r}-sets}
as the images $\bigl\{\phi(\alf):\alf\sle\{1,\ld,r\},\ \alf\neq\pus\bigr\}$
of polynomial mappings $\phi\col\cF(\{1,\ld,r\})\ra\Z^{l}$ 
of degree $\leq d$ and with $\phi(\pus)=0$,
and say that a set $E\in\Z^{l}$ is a {\it \VIPzr{d,r}-set\/}
if it has a nonempty intersection with every \VIPr{d,r}-set.
We will also say that a set $E\in\Z^{l}$ is a {\it \VIPzn-set\/}
if for any $d\in\N$, $E$ is an \VIPzr{d,r}-set for some $r\in\N$.
\rfr{P-Igpssub} now implies the following result:

\theorem{P-IpsVIPZ}{}{}
For any $D,d\in\N$ and $\eps>0$ there exists $r=r(D,d,\eps)\in\N$
such that for any $l\in\N$ 
and any open constant-free generalized polynomial mapping $\phi\col\Z^{l}\ra\R$ 
of total degree $\leq D$
the set $\bigl\{n\in\Z^{l}:\|\phi(n)\|<\eps\bigr\}$ is a \VIPzr{d,r}-set.
\endtheorem

Let $G$ be a nilpotent Lie group;
an {\it $l$-parameter polynomial sequence\/} in $G$ is a mapping $g\col\Z^{l}\ra G$
of the form $T_{1}^{p_{1}(n)}\cds T_{b}^{p_{b}(n)}$, $n\in\Z^{l}$,
where $T_{i}\in G$ and $p_{i}$ are polynomials $\Z^{l}\ra\Z$;
the {\it naive degree\/} of $g$ is defined as $\max_{i}\deg p_{i}$.%
\fnote{A more fundamental notion of {\it degree\/} of a polynomial sequence in a nilpotent group
can be defined as the number of ``differentiations'' 
which it takes in order to reduce the polynomial sequence to a constant.
For our purposes, however, the ``naive'' degree is quite sufficent.}
Using the fact that the coordinates of a point of a nilmanifold 
under the action of a polynomial sequence of niltranslations are generalized polynomials,
we obtain as a corollary of \rfr{P-IpsVIPZ} 
the following strengthening of the first part of \rfr{P-IPn}:

\theorem{P-IsubVIP}{}{(Cf.\nasp\ \rfr{P-nilRec} below.)}
Let $X$ be a nilmanifold with metric $\rho$
(compatible with the homogeneous space structure on $X$).
For any $a,d\in\N$ and $\eps>0$ there exists $r=r(a,d,\eps)\in\N$
such that for any $x_{0}\in X$,
any $l\in\N$, and any $l$-parameter polynomial sequence $g$ of niltranslations on $X$ 
of naive degree $\leq a$ and with $g(0)=\Id_{X}$,
the set $R_{U}(x_{0})=\bigl\{n\in\Z^{l}:\rho(g(n)x_{0},x_{0})<\eps\bigr\}$ is a \VIPzr{d,r}-set.
\endtheorem

We say that a dynamical system $(X,T)$ is {\it \VIPz-recurrent\/}
if for any $x_{0}\in X$ and any open neighborhood $U$ of $x_{0}$
the set of returns $R_{U}(x_{0})=\bigl\{n\in\Z:T^{n}x_{0}\in U\bigr\}$ is a \VIPz-set,
and is {\it\VIPzn-recurrent}
if for any $x_{0}\in X$ and any open neighborhood $U$ of $x_{0}$
the set $R_{U}(x_{0})$ is a \VIPzn-set.
The \VIPz-recurrence property turns out to be strictly stronger than that of the \IPz-recurrence:
there exist distal but not \VIPz-recurrent systems.%
\fnote{See \brfr{Ronnie}, Corollary~5.1, 
where it is shown that for any nonlinear polynomial $p\col\Z\ra\Z$
there exists an affine skew product transformation $T$
such that $\hbox{liminf}_{n}\hbox{dist}(T^{p(n)}0,0)>0$.}
As for the \VIPzn-recurrence, we get, as a corollary of \rfr{P-IsubVIP},
that, via \rfr{P-IPn}, \VIPzn-recurrence is equivalent to \IPzn-recurrence:
\theorem{P-VIPn}{}{}
Any pre-nilsystem is \VIPzn-recurrent,
and any \VIPzn-recurrent system is a disjoint union of pre-nilsystems.
\endtheorem
In \rfr{S-genpol} of the paper we prove (a more precise version of) 
Theorems~\rfrn{P-Igpssub} and \rfrn{P-IpsVIPZ} and deduce \rfr{P-IsubVIP} from them.
In \rfr{S-pronil} we obtain the  second statement of \rfr{P-IPn}.
\section{S-genpol}{Sets of visits of open bounded generalized polynomials with no constant term
to a neighborhood of zero.}

Let $A$ be a set and $(H,+)$ be an abelian group.
For $r\in\N$ we will denote by $\inr{r}$ the interval $\{1,\ld,r\}$ in $\N$.
We denote by $\cF(A)$ the set of finite subsets of $A$,
by $A^{(d)}$, $d\in\N$, the set of subsets of $A$ of cardinality $d$,
and by $A^{(\leq d)}$, $d\in\N$, 
the set of nonempty subsets of $A$ of cardinality $\leq d$,
$A^{(\leq d)}=\bigcup_{l=1}^{d}A^{(l)}$.

We start with discussing polynomial mappings on $\cF(A)$.
We say that a mapping $\phi\col\cF(A)\ra H$ is {\it linear\/}
if it satisfies the identity $\phi(\alf\cup\bet)=\phi(\alf)+\phi(\bet)$
whenever $\alf,\bet\in\cF(A)$ are disjoint,
and will denote the set of linear mappings $\cF(A)\ra H$ by $\Lin(A,H)$. 
A mapping $\phi\in\Lin(A,H)$ is uniquely defined by its values at singletons:
for any $\alf\in\cF(A)$, $\phi_{\alf}=\sum_{a\in\alf}\hphi(\{a\})$.
We will call the mapping $\hphi\col A\ra H$ defined by $\hphi(a)=\phi(\{a\})$
{\it the producing function\/} for $\phi$;
we then have $\phi(\alf)=\sum_{a\in\alf}\hphi(a)$, $\alf\in\cF(A)$.

For a mapping $\phi\col\cF(A)\ra H$ and $\bet\in\cF(A)$
we define {\it the $\bet$-derivative\/} $D_{\bet}\phi$ of $\phi$
by $D_{\bet}\phi(\alf)=\phi(\alf\cup\bet)-\phi(\alf)$, $\alf\in\cF(A\sm\bet)$.
We say that a mapping $\phi$ is {\it polynomial of degree $\leq d$}
if for any $d+1$ pairwise disjoint sets $\bet_{0},\ld,\bet_{d}\in\cF(A)$ 
one has $D_{\bet_{d}}\cds D_{\bet_{0}}\phi=0$.

We will denote by $\Pol_{d}(A,H)$
the group of polynomial mappings $\cF(A)\ra H$ of degree $\leq d$.
We will mainly deal with polynomial mappings ``having zero constant term'';
let us denote by $\Polz_{d}(A,H)$
the subgroup $\bigl\{\phi\in\Pol_{d}(A,H):\phi(\pus)=0\bigr\}$ of $\Pol_{d}(A,H)$.
Notice that $\Lin(A,H)=\Polz_{1}(A,H)$.

One can show (see \brfr{phj}, sections~8.3-8.5)
that any polynomial mapping $\phi\in\Polz_{d}(A,H)$
can be represented in the form $\phi(\alf)=\Phi(\alf^{d})$, $\alf\in\cF(A)$,
for some mapping $\Phi\in\Lin(A^{d},H)$,
so that 
\equn{f-cp}{
\phi(\alf)=\sum_{v\in\alf^{d}}\hPhi(v),\ \alf\in\cF(A), 
}
where $\hPhi\col A^{d}\ra H$ is the producing function for $\Phi$.
We will call $\hPhi$ {\it a \cbc-producing function for $\phi$}.

The \cbc-producing function for a polynomial mapping $\phi\in\Polz_{d}(A,H)$ 
is not canonically defined.
A more natural is {\it the \smp-producing function\/} for $\phi$, 
a function $\tPhi\col A^{(\leq d)}\ra H$ such that for any $\alf\in\cF(A)$, 
\equ{
\phi(\alf)=\sum_{u\in\alf^{(\leq d)}}\tPhi(u).
}
The \smp-producing function $\tPhi$ for $\phi$ is defined uniquely
(and provides a natural approach to the definition of polynomial mappings 
in the case $H$ is a commutative semigroup).
In terms of $\tPhi$,
$\phi$ is the sum of its {\it homogeneous components\/}, $\phi=\phi_{1}+\cds+\phi_{d}$,
where for each $i$, $\phi_{i}(\alf)=\sum_{\del\in\alf^{(d)}}\tPhi(\del)$.
To obtain the \smp-producing function $\tPhi$ for $\phi$
from a \cbc-producing function $\hPhi$
one simply sums up the values of $\hPhi$ at the elements of $A^{(d)}$ 
corresponding to the same element of $A^{(\leq d)}$:
for any $u\in A^{(\leq d)}$, 
\equn{f-cs}{
\tPhi(u)=\sum_{\sdup{v=(a_{1},\ld,a_{d})\in\alf^{d}}{\{a_{1},\ld,a_{d}\}=u}}\hPhi(v).
}

Let $B$ be a collection of pairwise disjoint finite subsets of $A$;
we will call $B$ {\it a \sbc{}\/} in $\PP{A}$;
if $|B|=s$ we will say that $B$ is {\it a \sbc{$s$}}.
Given a \sbc{} $B$ in $\PP{A}$,
we have an injection $\cF(B)\ra\cF(A)$ defined by $\gam\mapsto\bigcup\gam$,
and we will identify $\cF(B)$ with its image in $\cF(A)$.
Given a polynomial mapping $\phi\col\cF(A)\ra H$,
we call the polynomial mapping $\phi\rest{\cF(B)}$ {\it a subpolynomial\/} of $\phi$
corresponding to the \sbc{} $B$
and denote it by $\phi\rsub{B}$.
Any \sbc{} $B$ of a \sbc{} in $\PP{A}$
induces the \sbc{} $B'=\bigl\{\bigcup C:C\in B\bigr\}$ in $\PP{A}$;
abusing notation, we will denote the subpolynomial $\phi\rsub{B'}$ of $\phi$
by $\phi\rsub{B}$.

Let $\hPhi\col A^{d}\ra H$ be a \cbc-producing function
for a polynomial mapping $\phi\col\cF(A)\ra H$ of degree $\leq d$
and let $\Phi\in\Lin(A^{d},H)$ be the linear mapping produced by $\hPhi$.
Given a \sbc{$s$} $B=\{B_{1},\ld,B_{s}\}$ in $\PP{A}$,
one finds a \cbc-producing function for the subpolynomial $\phi\rsub{B}$ as follows.
For any $\bet\sle B$ we have
\equn{f-spc}{
\phi\rsub{B}(\bet)
\ &=\phi\Bigl(\bigcup_{C\in\bet}C\Bigr)
=\Phi\Bigl(\Bigl(\bigcup_{C\in\bet}C\Bigr)^{d}\Bigr)
=\sum_{v\in(\bigcup_{C\in\bet}C)^{d}}\hPhi(v)
\-\\
&=\sum_{C_{1},\ld,C_{d}\in\bet}\sum_{v\in C_{1}\times\cds\times C_{d}}\hPhi(v)
=\sum_{(C_{1},\ld,C_{d})\in\bet^{d}}\Phi(C_{1},\ld,C_{d});
\-}
thus, the mapping $\Phi\rest{B^{d}}$ is a \cbc-producing function for $\phi\rsub{B}$.

The following proposition establishes the \IPzr{r}-recurrence property
of polynomial mappings with values in the torus $\TT=\R/\Z$.
\proposition{P-subT}{}{(Cf.\nasp\ \brfr{B-ultra}, Theorem~7.7)}
For any $k,d\in\N$ and $\eps>0$ there exists $r=r(k,d,\eps)\in\N$ such that
for any $\phi_{1},\ld,\phi_{k}\in\Polz_{d}(\inr{r},\TT)$
there exists a nonempty $\alf\in\cF(\inr{r})$
such that $\dist(\phi_{i}(\alf),0)<\eps$ for all $i\in\{1,\ld,k\}$
(where ``$\dist$'' is the distance on $\TT$).
\endproposition
\proof{}
Put $c=\lceil1/\eps\rceil$
and partition the torus $\TT$ into $c$ intervals of length $\leq 1/\eps$.
By the Polynomial Hales-Jewett theorem (see \brfr{phj}, Theorem~0.10),
there exists $r\in\N$ such that for any partition of $\cF(\inr{r}^{d}\times[k])$ into $c$ subsets
there exist $\gam\sln\inr{r}^{d}\times[k]$ and a nonempty $\alf\sle\inr{r}$
such that $\gam\cap(\alf^{d}\times[k])=\pus$ 
and the sets $\gam,\gam\cup(\alf^{d}\times\{1\}),\dsc\ld,\gam\cup(\alf^{d}\times\{k\})$ 
belong to the same element of the partition.
Let $\phi_{1},\ld,\phi_{k}\in\Polz_{d}(\inr{r},\TT)$. 
For each $i$ let $\hPhi_{i}\col\inr{r}^{d}\ra\TT$ be a \cbc-producing function for $\phi_{i}$.
Define a mapping $\hPhi\col\inr{r}^{d}\times[k]\ra\TT^{k}$
by $\hPhi(v,i)=\hPhi_{i}(v)$, $v\in\inr{r}^{d}$, $i\in[k]$, 
and let $\Phi\in\Lin(\inr{r}^{d}\times[k],\TT)$ be the linear mapping produced by $\hPhi$.
Then, via $\Phi$, the partition of $\TT$ 
defines a partition of $\cF(\inr{r}^{d}\times[k])$ into $c$ subsets.
Applying the Polynomial Hales-Jewett theorem,
we can find $\gam\sln\inr{r}^{d}\times[k]$ and a nonempty $\alf\sle\inr{r}$ 
such that $\gam\cap(\alf^{d}\times[k])=\pus$ 
and the sets $\gam,\gam\cup(\alf^{d}\times\{1\}),\ld,\gam\cup(\alf^{d}\times\{k\})$
belong to the same element of the partition;
then for any $i$,
$\Phi(\gam)$ and $\Phi(\gam\cup(\alf^{d}\times\{i\}))$ belong to the same partition of $\TT$,
and so, $\dist\bigl(\Phi(\gam),\Phi(\gam\cup(\alf^{d}\times\{i\}))\bigr)<\eps$.
Since $\Phi(\gam\cup(\alf^{d}\times\{i\}))=\Phi(\gam)+\Phi(\alf^{d}\times\{i\})
=\Phi(\gam)+\phi_{i}(\alf)$,
this implies that $\dist(0,\phi_{i}(\alf))<\eps$.%
\endproof
Recall that by $\lip x\rip$ we denote the distance from $x\in\R$ to $\Z$.
We may then reformulate \rfr{P-subT} as follows:
\corollary{P-sub}{}{}
For any $k,d\in\N$ and $\eps>0$ there exists $r=r(k,d,\eps)\in\N$ such that
for any $\phi_{1},\ld,\phi_{k}\in\Polz_{d}(\inr{r},\R)$
there exists a nonempty $\alf\in\cF(\inr{r})$ 
such that $\lip\phi_{i}(\alf)\rip<\eps$ for all $i\in\{1,\ld,k\}$.
\endcorollary
Next we show that if $r$ is large enough,
any polynomial mapping $\phi\in\Polz_{d}(\inr{r},\TT)$ has a subpolynomial 
whose \cbc-producing function is arbitrarily small:
\proposition{P-subcmT}{}{}
For any $d,s\in\N$ and $\eps>0$ there exists $r\in\N$ such that
for any $\phi\in\Polz_{d}(\inr{r},\TT)$
there exists a \sbc{$s$} $B$ in $\PP{\inr{r}}$
such that a \cbc-producing function $\hPhi_{B}$ for $\phi\rsub{B}$
satisfies $\dist(\hPhi_{B},0)<\eps$.
\endproposition
\proof{}
Take $r_{0}=r(s^{d},d,\eps)$ as in \rfr{P-sub},
and put $A=\inr{s}\times\inr{r_{0}}$ 
(and $r=|A|=sr_{0}$).
Let $\phi\in\Polz_{d}(\inr{r},\TT)$.
Let $\hPhi\col\inr{r}\ra\TT$ be a \cbc-producing function for $\phi$
and let $\Phi\in\Lin(\inr{r}^{d},\TT)$ be the linear mapping produced by $\hPhi$.
For each $i=(i_{1},\ld,i_{d})\in\inr{s}^{d}$
define a polynomial mapping $\phi_{I}\in\Polz_{d}(\inr{r_{0}},\TT)$
by $\phi_{I}(\alf)=\Phi((\{i_{1}\}\times\alf)\times\cds\times(\{i_{d}\times\alf)\})$.
By \rfr{P-sub} there exists $\alf\sle\inr{r_{0}}$
such that $\dist(\phi_{I}(\alf),0)<\eps$ for all $I\in\inr{s}^{d}$.
Take the \sbc{$s$} $B=\bigl\{\{i\}\times\alf:i\in\inr{s}\bigr\}$ in $\PP{A}$.
By the choice of $\alf$,
for any $w\in B^{d}$ we have $\dist(\Phi(w),0)<\eps$.
Since, by \frfr{f-spc}, $\Phi\rest{B^{d}}$ is a \cbc-producing function for $\phi\rsub{B}$,
we are done.
\endproof
Replacing in \rfr{P-subcmT} $\eps$ by $\eps/s^{d}$, we obtain:
\corollary{P-subtmT}{}{}
For any $d,s\in\N$ and $\eps>0$ there exists $r\in\N$ such that
for any $\phi\in\Polz_{d}(\inr{r},\TT)$
there exists a \sbc{$s$} $B$ in $\PP{\inr{r}}$
such that $\dist(\phi\rsub{B},0)<\eps$.
\endcorollary
By formula \frfr{f-cs},
any value of the \smp-producing function for $\phi\in\Polz_{d}(A,\R)$ 
is a sum of less than $d^{d}$ values of the \cbc-producing function for $\phi$.
Hence, \rfr{P-subcmT} implies the following corollary:
\proposition{P-subsmT}{}{}
For any $d,s\in\N$ and $\eps>0$ there exists $r\in\N$ such that
for any $\phi\in\Polz_{d}(\inr{r},\TT)$
there exists a \sbc{$s$} $B$ in $\PP{\inr{r}}$
such that the \smp-producing function $\tPhi_{B}$ for $\phi\rsub{B}$
satisfies $\dist(\tPhi_{B},0)<\eps$.
\endproposition
In terms of polynomial mappings with values in $\R$, 
\rfr{P-subsmT} takes the following form:
\corollary{P-subsm}{}{}
For any $d,s\in\N$ and $\eps>0$ there exists $r\in\N$ such that
for any $\phi\in\Polz_{d}(\inr{r},\R)$
there exists a \sbc{$s$} $B$ in $\PP{\inr{r}}$
such that the \smp-producing function $\tPhi_{B}$ for $\phi\rsub{B}$
satisfies $\lip\tPhi_{B}\rip<\eps$.
\endproposition
Now let $\phi\in\Polz_{d}(\inr{r},\R)$ be a polynomial mapping
whose \smp-producing function $\tPhi$ satisfies $\lip\tPhi\rip<1/r^{d}$.
We will denote by $[x]$ the integer and by $\{x\}$ the fractional parts of $x\in\R$.
If $x\in\R$ satisfies $\lip x\rip<\eps$, then either $\{x\}<\eps$ or $\{x\}>1-\eps$.
If $x_{1},\ld,x_{n}\in\R$ satisfy $\{x_{i}\}<1/n$, $i=1,\ld,n$,
then $\bigl[\sum_{i=1}^{n}x_{i}\bigr]=\sum_{i=1}^{n}[x_{i}]$.
Thus, if $\tPhi$ satisfies $\{\tPhi\}<1/r^{d}$,
then for any $\alf\sle\inr{r}$, 
\equ{
[\phi(\alf)]=\Bigl[\sum_{u\in\alf^{(\leq d)}}\tPhi(u)\Bigr]
=\sum_{u\in\alf^{(\leq d)}}[\tPhi(u)]
}
and so, $[\phi]$ is also a polynomial mapping, $[\phi]\in\Polz_{d}(\inr{r},\Z)$,
with the \smp-producing function $[\tPhi]$.

For any $x\in\R\sm\Z$, $[x]=-[-x]-1$ and $\{-x\}=1-\{x\}$,
so, if $x_{1},\ld,x_{n}\in\R$ satisfy $\{x_{i}\}>1-1/n$, $i=1,\ld,n$,
then 
$$\txt
\bigl[\sum_{i=1}^{n}x_{i}\bigr]
=-\bigl[-\sum_{i=1}^{n}x_{i}\bigr]-1
=-\bigl[\sum_{i=1}^{n}(-x_{i})\bigr]-1
=-\sum_{i=1}^{n}[-x_{i}]-1
=\sum_{i=1}^{n}(-[-x_{i}])-1.
$$
Applying this to $\tPhi$,
we see that if $\tPhi$ satisfies $\{\tPhi\}>1-1/r^{d}$,
then for any $\alf\sle\inr{r}$, 
\equ{
[\phi(\alf)]=\Bigl[\sum_{u\in\alf^{(\leq d)}}\tPhi(u)\Bigr]
=\sum_{u\in\alf^{(\leq d)}}\Bigl(-[-\tPhi(u)]\Bigr)-1.
}
So, $[\phi]+1$ is a polynomial mapping, $[\phi]+1\in\Polz_{d}(\inr{r},\Z)$,
with the \smp-producing function $-[-\tPhi]$.

In the general case, when $\lip\tPhi\rip<1/r^{d}$,
we may have neither $\{\tPhi\}<1/r^{d}$ nor $\{\tPhi\}>1-1/r^{d}$.
However, if $\phi$ is a homogeneous polynomial of degree $l\leq d$
(which means that $\phi(\alf)=\sum_{u\in\alf^{(l)}}\tPhi(u)$),
then, given $s\in\N$, if $r$ is large enough,
by the classical Ramsey theorem we can choose an $s$-element subset $B$ of $\inr{r}$
such that either $\{\tPhi(u)\}<1/r^{d}$ for all $u\in B^{(d)}$
or $\{\tPhi(u)\}>1-1/r^{d}$ for all $u\in B^{(d)}$.
Identifying $B$ with the ``singleton \sbc'' $\{\{b\}:b\in B\}$ in $\PP{\inr{r}}$,
we will therefore have $[\phi\rsub{B}]\in\Polz_{d}(B,\Z)+e$
with $e\in\{0,-1\}$.

For a general $\phi\in\Polz_{d}(\inr{r},\R)$,
applying this argument to all homogeneous components of $\phi$ and using a diagonal process,
we arrive at the following lemma:
\lemma{P-skob}{}{}
For any $d,s\in\N$ there exists $r\in\N$
such that for any $\phi\in\Polz_{d}(\inr{r},\R)$ 
whose \smp-producing function $\tPhi$ satisfies $\lip\tPhi\rip<1/r^{d}$
there exists a (singleton) \sbc{} $B$ in $\PP{\inr{r}}$
such that $[\phi\rsub{B}]\in\Polz_{d}(B,\Z)+e$ with $e\in\{0,-1,\ld,-d\}$.
\endlemma
Combining \rfr{P-skob} with \rfr{P-subsm} we obtain:
\theorem{P-skop}{}{}
For any $d,s\in\N$ there exists $r\in\N$ 
such that for any $\phi\in\Polz_{d}(\inr{r},\R)$ there exists a \sbc{$s$} $B$ in $\PP{\inr{r}}$
such that $[\phi]\in\Polz_{d}(B,\Z)+e$ with $e\in\{0,-1,\ld,-d\}$.
\endtheorem
Using induction on $k$, one can extend \rfr{P-skop} to the case of $k$ polynomials:
\theorem{P-skopk}{}{}
For any $k,d_{1},\ld,d_{k},s\in\N$ there exists $r=r(k,(d_{1},\ld,d_{k}),s)\in\N$ 
such that for any $\phi_{i}\in\Polz_{d_{i}}(\inr{r},\R)$, $i=1,\ld,k$,
there exists a \sbc{$s$} $B$ in $\PP{\inr{r}}$
such that for every $i\in\{1,\ld,k\}$,
$[\phi_{i}]\in\Polz_{d_{i}}(B,\Z)+e_{i}$, with $e_{i}\in\{0,-1,\ld,-d\}$.
\endtheorem
\vbreak{-9000}{4mm}
{\it A generalized polynomial\/} is a function obtained from conventional polynomials
using the operations of taking the integer part, addition, and multiplication.
We say that a generalized polynomial $\phi$ is {\it constant free\/}
if all polynomials involved in the expression of $\phi$
have zero constant term.
(More precisely, a generalized polynomial is constant free if it has {\it a representation\/}
in which all polynomials have zero constant term.
A similar convention applies to all the definitions below.)
We say that a polynomial $\phi$ is {\it open\/}
if it is contained in the ideal, in the ring of constant free generalized polynomials,
generated by the ordinary polynomials.
This is equivalent to saying that $\phi$
(or rather a representation of $\phi$) has no summand 
that is a product of ``closed'' generalized polynomials $[\phi_{i}]$.
Any open constant-free generalized polynomial is representable in the form
\equn{f-gp}{
\phi=\sum_{j=1}^{m}[\phi_{j,1}]\cds[\phi_{j,l_{j}}]\phi_{j,0}
}
where for every $j$,
$\phi_{j,1},\ld,\phi_{j,l_{j}}$ are open constant-free generalized polynomials
and $\phi_{j,0}$ are conventional polynomials with zero constant term.

We now introduce the notions of height, width, and degree 
for (a representation of) a generalized polynomial $\phi$:\\
{\it The height\/} $h(\phi)$ of $\phi$ is the maximum length 
of sequences of nested brackets in $\phi$:
we put $h(\phi)=0$ if $\phi$ is a conventional polynomial
and we say that $h(\phi)\leq h$ if $\phi$ has a representation \frfr{f-gp}
where for all $j$ and all $t\geq 1$, $h(\phi_{j,t})\leq h-1$.\\
{\it The width\/} $w(\phi)$ is the maximum number of components in $\phi$ itself
and in all its components:
we put $w(\phi)=1$ if $\phi$ is a conventional polynomial
and we say that $w(\phi)\leq w$ if $\phi$ has a representation \frfr{f-gp}
where $w(\phi_{j,t})\leq w$ for all $j$ and all $t\geq 1$ 
and also $\sum_{j=1}^{m}(l_{j}+1)\leq w$.\\
{\it The degree\/} $d(\phi)$ of $\phi$ is defined as usual 
under the assumption that $\deg[\phi]=\deg\phi$:
we say that $d(\phi)\leq d$ if $\phi$ has a representation \frfr{f-gp}
with $\max_{j=1}^{m}\bigl(\sum_{t=0}^{l_{j}}\deg\phi_{j,t}\bigr)\leq d$.\\
(For example, for $\phi(x)=[[x^{2}+1]x][x^{3}+2x]x+[x^{2}](x+1)+x^{3}$ 
we have $h(\phi)=2$, $w(\phi)=6$, and $d(\phi)=7$.)

We extend the above definitions to generalized polynomial mappings with domain $\cF(A)$,
and will denote by $\GPolzo_{d,h,w}(A,H)$ the set (the algebra) 
of open constant-free generalized polynomial mappings $\phi\col\cF(A)\ra H$, where $H=\R$ or $\Z$,
with $d(\phi)\leq d$, $h(\phi)\leq h$, and $w(\phi)\leq w$.
Given $\phi\in\GPolzo_{d,h,w}(A,H)$ and a \sbc{} $B$ in $\PP{A}$,
we define the generalized polynomial mapping $\phi\rsub{B}\in\GPolzo_{d,h,w}(B,H)$
as the restriction of $\phi$ to the set $\cF(B)$ considered as a subset of $\cF(A)$.

The following theorem says that generalized polynomial mappings
turn into ordinary polynomial mappings
after being restricted to a suitable \sbc{} in their domain:
\theorem{P-gps}{}{}
For any $k,d_{1},\ld,d_{k},h,w,s\in\N$ 
there exists $r=r(k,(d_{1},\ld,d_{k}),h,w,s)\dsc\in\N$ such that
for any $\phi_{i}\in\GPolzo_{d_{i},h,w}(\inr{r},\R)$, $i=1,\ld,k$,
there exists a \sbc{$s$} $B$ in $\PP{\inr{r}}$
such that $\phi_{i}\rsub{B}\in\Polz_{d_{i}}(B,\R)$, $i=1,\ld,k$.
\endtheorem
\proof{}
We will use induction on $h$; when $h=0$ the statement is trivial.
Take $r_{0}$ to be the maximum of the integers $r(l,(b_{1},\ld,b_{l}),s)$ in \rfr{P-skopk}
over all integers $l\leq kw$ 
and all $l$-tuples $(b_{1},\ld,b_{l})$ of nonnegative integers 
with $\sum_{j=1}^{l}b_{j}\leq w\sum_{i=1}^{k}d_{i}$.
By induction on $h$, 
let $r$ be the maximum of the integers $r(l,(d_{1},\ld,d_{l}),h-1,w,r_{0})$ 
in the assertion of \rfr{P-gps}
over all integers $l\leq kw$
and all $l$-tuples $(b_{1},\ld,b_{l})$ of nonnegative integers 
with $\sum_{j=1}^{l}b_{j}\leq w\sum_{i=1}^{k}d_{i}$.
Let $\phi_{i}\in\GPolzo_{d_{i},h,w}(\inr{r},\R)$, $i=1,\ld,k$.
For each $i$ reprsent $\phi_{i}$ in the form
$$
\phi_{i}=\sum_{j=1}^{m_{i}}[\phi_{i,j,1}]\cds[\phi_{i,j,l_{i,j}}]\phi_{i,j,0},
$$
where for every $i,j$ we have $\phi_{i,j,0}\in\Polz_{d_{i,j,0}}(\inr{r},\R)$ 
and for every $t\geq 1$ we have $\phi_{i,j,t}\in\GPolzo_{d_{i,j,t},h-1,w}(\inr{r},\R)$
with 
$$\hbox{
$\sum_{j=1}^{m_{i}}(l_{i,j}+1)\leq w$ for all $i$
and $\sum_{t=0}^{l_{i,j}}d_{i,j,t}\leq d_{i}$ for all $i,j$,
}$$
so that 
$$\hbox{
$\sum_{i=1}^{k}\sum_{j=1}^{m_{i}}(l_{i,j}+1)\leq kw$ 
and $\sum_{i=1}^{k}\sum_{j=1}^{m_{i}}\sum_{t=0}^{l_{i,j}}d_{i,j,t}\leq w\sum_{i=1}^{k}d_{i}$.
}$$
By the choice of $r$
there exists a \sbc{$r_{0}$} $B_{0}\sln\cF(\inr{r})$
such that $\phi_{i,j,t}\rsub{B_{0}}\in\Polz_{d_{i,j,t}}(B_{0},\R)$ for all $i,j,t$.
Then by the choice of $r_{0}$ there exists a \sbc{$s$} $B$ in $\PP{B_{0}}$
such that for all $i,j,t$,
$\bigl[\phi_{i,j,t}\rsub{B}\bigr]\in\Polz_{d_{i,j,t}}(B,\Z)$.
Hence for every $i$,
$$
\phi_{i}\rsub{B}=\sum_{j=1}^{m_{i}}\bigl[\phi_{i,j,1}\rsub{B}\bigr]\cds
\bigl[\phi_{i,j,l_{i,j}}\rsub{B}\bigr]\phi_{i,j}\rsub{B}
\in\Polz_{d_{i}}(B,\R).
$$
\frgdsp\endproof
Combining \rfr{P-gps} and \rfr{P-sub}, we obtain:

\theorem{P-gpssub}{}{}
For any $k,d,h,w\in\N$
there exists $r=r(k,d,h,w)\in\N$ such that
for any $\phi_{1},\ld,\phi_{k}\in\GPolzo_{d,h,w}(\inr{r},\R)$
there exists a nonempty $\alf\in\cF(\inr{r})$
such that $\lip\phi_{i}(\alf)\rip<\eps$, $i=1,\ld,k$.
\endtheorem

Let $X=G/\Gam$ be a $k$-dimensional compact nilmanifold; 
we may and will assume that $X$ is connected.
(Any nilmanifold is a subnilmanifold of a connected one.)
Let $\rho$ be a metric on $X$ 
(induced by a metric on $G$ compatible with the Lie group structure thereon).
Fix a point $x_{0}\in X$,
and let $\tau=(\tau_{1},\ld,\tau_{k})\col X\ra[0,1)^{k}$ 
be Maltsev's coordinates on $X$ centered at $x_{0}$.
The inverse mapping $\tau^{-1}$ is continuous,
and the distance $\rho(x,x_{0})$ from $x\in X$ to $x_{0}$
is continuous with respect to the distance 
from $\tau(x)$ to the set of vertices $\{0,1\}^{k}$ of the cube $[0,1]^{k}$.
(See, for example, \brfr{sko}, Section~1.5.)

Let $g$ be an {\it ($l$-parameter) polynomial sequence\/} in $G$,
that is, a mapping $g\col\Z^{l}\ra G$ 
of the form $g(n)=T_{1}^{p_{1}(n)}\cds T_{b}^{p_{b}(n)}$, $n\in\Z^{l}$,
where $T_{1},\ld,T_{b}\in G$,
$p_{1},\ld,p_{b}$ are polynomials $\Z^{l}\ra\Z$;
we define $\ndeg g$, {\it the naive degree of $g$\/}, as $\max_{i=1}^{b}\deg p_{i}$.
Then for each $i=1,\ld,k$, 
the sequences $\psi_{i}(n)=\tau_{i}(g(n)x_{0})$, $n\in\Z^{c}$, 
of coordinates of $x_{0}$ under the action of $g$
are open $[0,1)$-valued generalized polynomials,
with parameters depending only on $X$ and $\ndeg g$
(see \brfr{sko}, Theorem~A and Theorem~A$^{**}$),
and if $g(0)=1_{G}$, these polynomials can be assumed to be constant-free.
For any polynomial mapping $\phi\in\Pol_{d}(\inr{r}),\Z^{c})$,
the composition mappings $\psi_{i}\comp\phi\col\cF(\inr{r})\ra[0,1)$, $i=1,\ld,k$,
are open constant-free generalized polynomial mappings,
with parameters only depending on $X$, $d$, and $\ndeg g$.
From \rfr{P-gpssub} we now obtain the following result:

\theorem{P-nilRec}{}{}
Let $X=G/\Gam$ be a nilmanifold with metric $\rho$.
For any $a,d\in\N$ and $\eps>0$ there exists $r=r(a,d,\eps)\in\N$ such that
for any $l$, any $l$-parameter polynomial sequence $g$ in $G$
with $\ndeg g\leq a$ and $g(0)=1_{G}$,
any $x_{0}\in X$, and any $\phi\in\Polz_{d}(\inr{r},\Z^{l})$
there exists a nonempty $\alf\in\cF(\inr{r})$ 
such that $\rho\bigl(g(\phi(\alf))x_{0},x_{0}\bigr)<\eps$.
\endtheorem
\remark{}
\rfr{P-nilRec} easily extends to {\it generalized polynomial sequences\/} in nilpotent groups,
that is, to sequences of the form $g(n)=T_{1}^{p_{1}(n)}\cds T_{b}^{p_{b}(n)}$
where $p_{i}$ are generalized polynomials $\Z^{l}\ra\Z$.
\endremark
\section{S-pronil}{\IPzn-recurrence implies approximability by nilsystems}

In this section we prove the second statement of \rfr{P-IPn}.
Let $(X,\ro)$ be a compact metric space,
$T$ be a self homeomorphism of $X$,
and assume that $(X,T)$ is \IPzn-recurrent.
Then, in particular, $(X,T)$ is \IPz-recurrent,
so by \rfr{P-IP}, $(X,T)$ is distal,
and thus is a disjoint union of minimal subsystems 
(see \brfr{F-book}, corollary to Theorem~8.7).
Hence, we may assume that $(X,T)$ is minimal.

Now, by the way of contradiction, 
assume that a minimal system $(X,T)$ is not a pre-nilsystem,
that is, not an inverse limit of nilsystems;
our goal is to show that there exists a point $x\in X$ and $\eps>0$
such that for every $r\in\N$ there exists a linear mapping $\phi\in\Lin(\inr{r},\Z)$
such that $\ro(T^{\phi(\alf))}x,x)>\eps$ for every nonempty $\alf\sle\inr{r}$.

We will use the following result (\brfr{HKM} Theorem~1.3 and Corollary~4.2):
for any $r$, the maximal $r$-step pro-nilfactor of $(X,T)$
is defined by a closed $T$-invariant equivalence relation $\RP{r}\sle X^{2}$
(called {\it the regionally proximal relation of order $r$}),
with $(x_{0},y_{0})\in\RP{r}$ if and only if for any $\del>0$ there exists a point $x\in X$ 
and a mapping $\phi\in\Lin(\inr{r},\Z)$
such that 
\equn{f-Pr}{
\hbox{$\ro(x,x_{0})<\del$ and $\ro(T^{\phi(\alf)}x,y_{0})<\del$ 
for all nonempty $\alf\sle\inr{r}$.}
}
Our assumption that $(X,T)$ is not a pre-nilsystem 
is equivalent to the assumption that $\bigcap_{r=1}^{\infty}\RP{r}\neq\Del$,
where $\Del$ is the diagonal of $X^{2}$.
Fix $(x_{0},y_{0})\in\bigcap_{r=1}^{\infty}\RP{r}$ with $x_{0}\neq y_{0}$.
Let $\eps=\inf_{n\in\Z}\ro(T^{n}x_{0},T^{n}y_{0})$;
since $(X,T)$ is distal, we have $\eps>0$.
Since $(X,T)$ is minimal,
the orbit $\{T^{n}x_{0}\}_{n\in\Z}$ of $x_{0}$ is dense in $X$.
Let $r\in\N$ and let $U\sle X$ be an open set.
Choose $n\in\Z$ such that $T^{n}x_{0}\in U$
and choose $\del>0$ such that $\ro(T^{n}x,T^{n}y)<\eps/3$ whenever $\ro(x,y)<\del$.
Find $x\in X$ such that \frfr{f-Pr} holds and $T^{n}x\in U$.
Then $\ro(T^{n}x,T^{n}x_{0})<\eps/3$ and 
$\ro(T^{\phi(\alf)}T^{n}x,T^{n}y_{0})<\eps/3$ for all nonempty $\alf\sle\inr{r}$,
and since $\ro(T^{n}x_{0},T^{n}y_{0})\geq\eps$,
we have that $\ro(T^{\phi(\alf)}T^{n}x,T^{n}x)>\eps/3$ for all nonempty $\alf\sle\inr{r}$.
This proves that for any $r\in\N$ the open set 
\equ{\txt
R_{r}=\bigl\{x\in X:\hbox{there exists $\phi\in\Lin(\inr{r},\Z)$ 
such that $\ro(T^{\phi(\alf)}x,x)>\eps/3$}
\kern2cm\\\-\txt
\hbox{for all nonempty $\alf\in\inr{r}$}\Bigr\}
}
is dense in $X$.
By Baire category theorem $\bigcap_{r=1}^{\infty}R_{r}$ is nonempty,
which gives us what we wanted --
a point $x\in X$
such that for every $r\in\N$ there exists a mapping $\phi\in\Lin(\inr{r},\Z)$
such that $\ro(T^{\phi(\alf))}x,x)>\eps/3$ for every nonempty $\alf\sle\inr{r}$.
\bibliography{}
\biblleft=13mm
\bibook AGH/AGH
a:L. Auslander, L. Green and F. Hahn
t:Flows on homogeneous spaces
i:Annals of Math. Studies, vol.~53,
Princeton Univ. Press, 1963
*
\bibart AGN/AGN
a:L. Auslander, G. Greschonig and A. Nagar
t:Reflections on equicontinuity
j:Proc. of AMS
n:142
y:2014
p:\no{9}, 3129-3137
*
\bibart B-ultra/B
a:V. Bergelson
t:Ultrafilters, IP sets, Dynamics, and Combinatorial Number Theory
j:Contemp. Math.
n:530
y:2010
p:23-47
*
\bibart phj/BL1
a:V. Bergelson and A. Leibman
t:Set-polynomials and polynomial extension of the Hales-Jewett Theorem
j:Ann. of Math.
n:150
y:1999
p:\no{1}, 33-75 
*
\bibart sko/BL2
a:V. Bergelson and A. Leibman
t:Distribution of values of bounded generalized polynomials
j:Acta Math.
n:198
y:2007
p:155-230
*
\bibarta bic/BL3
a:V. Bergelson and A. Leibman
t:Sets of large values of correlation functions for polynomial cubic configurations
j:Ergodic Th. and Dynam. Sys.
*
\bibook F-book/F
a:H. Furstenberg
t:Recurrence in Ergodic Theory and Combinatorial Number Theory
i:Princeton Univ. Press, 1981
*
\bibart HKM/HKM
a:B. Host, B. Kra, and A. Maass
t:Nilsequences and a structure theorem for topological dynamical systems
j:Adv. in Math.
n:224
y:2010
p:\no{1}, 103-129
*
\bibart Keyns1/Ke1
a:H. B. Keynes
t:Topological dynamics in coset transformation groups
j:Bull. Amer. Math. Soc.
n:72
y:1966
p:1033-1035
*
\bibart Keyns2/Ke2
a:H. B. Keynes
t:A study of the proximal relation in coset transformation groups 
j:Trans. Amer. Math. Soc.
n:128
y:1967
p:389-402
*
\bibart Ronnie/P
a:R. Pavlov
t:Some counterexamples in topological dynamics
j:Ergod. Th. and Dynam. Sys.
n:28
y:2008
p:\no{4}, 1291-1322
*
\bibart Veech/V
a:W. A. Veech
t:The equicontinuous structure relation for minimal Abelian transformation groups
j:Amer. J. Math.
n:90
y:1068
p:723-732
*
\endbibliography
\finish
\bye